\documentclass{article}

\newtheorem{theorem}{Theorem}[section]
\newtheorem{corollary}[theorem]{Corollary}

\newtheorem{lemma}[theorem]{Lemma}
\newtheorem{proposition}[theorem]{Proposition}
\newtheorem{definition}[theorem]{Definition}

 \evensidemargin\oddsidemargin
\begin{document}

\title{Correspondence between two antimatroid algorithmic characterizations }
\author{Yulia Kempner and Vadim E. Levit \\
Department of Computer Science\\
Holon Academic Institute of Technology\\
52 Golomb Str., P.O. Box 305\\
Holon 58102, ISRAEL\\
yuliak,levitv@hait.ac.il}
\date{}
\maketitle

\begin{abstract}
The basic distinction between already known algorithmic characterizations of
matroids and antimatroids is in the fact that for antimatroids the ordering
of elements is of great importance.

While antimatroids can also be characterized as set systems, the question
whether there is an algorithmic description of antimatroids in terms of sets
and set functions was open for some period of time.

This article provides a selective look at classical material on algorithmic
characterization of antimatroids, i.e., the ordered version, and a new
unordered version. Moreover we empathize formally the correspondence between
these two versions.\newline

\textbf{keywords:}\textit{\ antimatroid, greedoid, chain algorithm, greedy
algorithm, monotone linkage function.}
\end{abstract}

\section{Introduction}

In this paper we compare two algorithmic characterization of antimatroids.
There are many equivalent axiomatizations of antimatroids, that may be
separated into two categories: antimatroids defined as set systems and
antimatroids defined as languages. Boyd and Faigle \cite{BF} introduced an
algorithmic characterization of antimatroids based on the second antimatroid
definition - as a formal language. Another characterization of antimatroids,
that uses their definition as a set systems, is considered in this paper.
This approach is based on the optimization of set functions defined as the
minimum value of linkages between a set and elements from the set
complement. The correspondence between two these approaches is established.

Section 2 gives some basic information about antimatroids as set systems. In
Section 3 a set system generated by an isotone operator is introduced, and
its equivalence to an antimatroid is proved. In Section 4 monotone linkage
functions are considered. The optimization of the functions defined as the
minimum of the monotone linkage functions extends to antimatroids, and the
polynomial algorithm that finds an optimal set is constructed. In Section 5
an algorithmic characterization of truncated antimatroids in terms of the
monotone linkage functions are considered. In Section 6 the results of Boyd
and Faigle are presented and connection between their approach and the
approach based on monotone linkage functions is established.

\section{Preliminaries}

Let $E$ be a finite set. A \textit{set system }over $E$\textit{\ }is a pair $%
(E,\mathcal{F})$ where $\mathcal{F}\subseteq 2^{E}$ is a family of subsets
of $E$, called \textit{feasible }sets. We will use $X\cup x$ for $X\cup
\{x\} $ and $X-x$ for $X-\{x\}$.

\begin{definition}
A non-empty set system $(E,\mathcal{F})$ is an \textit{antimatroid }if

$(A1)$ for each non-empty $X\in \mathcal{F}$ there is an $x\in X$ such that $%
X-x\in \mathcal{F}$

$(A2)$ for all $X,Y\in \mathcal{F}$, and $X\not{\subseteq }Y$, there exist
an $x\in X-Y$ such that $Y\cup x\in \mathcal{F}$.
\end{definition}

\smallskip Any set system satisfying $(A1)$ is called \textit{accessible.}

\begin{definition}
\smallskip A set system $(E,\mathcal{F})$ has the \textit{interval property
without upper bounds }if for all $X,Y\in \mathcal{F}$ with $X\subseteq Y$
and for all $x\in E-Y$ , $X\cup x\in \mathcal{F}$ implies $Y\cup x\in 
\mathcal{F}$.
\end{definition}

\smallskip There are some different antimatroid definitions, for the sake of
completeness we will prove the following proposition:

\begin{proposition}
\label{P-1}\cite{BZ}\cite{Greedoids}For an accessible set system $(E,%
\mathcal{F})$ the following statements are equivalent:

$(i)$ $(E,\mathcal{F})$ is an antimatroid

$(ii)$ $\mathcal{F}$ is closed under union

$(iii)$ $(E,\mathcal{F})$ satisfies the interval property without upper
bounds
\end{proposition}

\setlength {\parindent}{0.0cm}\textbf{Proof.} $(i)\Rightarrow (ii)$ Let $%
X,Y\in \mathcal{F}$ and $Y\not{\subseteq }X$. Repeated application of $(A2)$
yields a set $X\cup (Y-X)\in \mathcal{F}$, i.e. $X\cup Y\in \mathcal{F}$. 
\setlength
{\parindent}{3.45ex}

$(ii)\Rightarrow (iii)$ Let $X,Y\in \mathcal{F}$ and $X\subseteq Y$ and $%
x\in E-Y$ and $X\cup x\in \mathcal{F}$, then $(X\cup x)\cup Y$ $\in \mathcal{%
F}$, i.e. $Y\cup x\in \mathcal{F}$.

$(iii)\Rightarrow (i)$ Let $X,Y\in \mathcal{F}$ and $X\not{\subseteq }Y$ .
Accessibility means that we can find a sequence $x_{1}x_{2}...x_{k}$ and
corresponding sequence $\emptyset =X_{0}\subset X_{1}\subset ...\subset
X_{k}=X$ where $X_{i}=X_{i-1}\cup x_{i}$ , and $X_{i}\in \mathcal{F}$ for $%
0\leq i\leq k$. Let $j$ be the least integer for which $X_{j}\not{\subseteq
}Y$. Then $X_{j-1}\subseteq Y$, $x_{j}\notin Y$ and $X_{j-1}\cup x_{j}\in 
\mathcal{F}$ , that implies $Y\cup x_{j}\in \mathcal{F}$. Hence $(E,\mathcal{%
F})$ is an antimatroid. \rule{2mm}{2mm}

\smallskip A maximal feasible subset of set $X\subseteq E$ is called a 
\textit{basis of }$X$, and will be\textit{\ }denoted by \-$\mathcal{B}_{X}$.
\smallskip Clearly, by $(ii)$, there is only one basis for each set.

\smallskip For a set $X\in \mathcal{F}$, let $\Gamma (X)=\{x\in E-X:X\cup
x\in \mathcal{F}\}$ be the set of \textit{feasible continuations }of $X$ 
\cite{Greedoids}.

We will say that $\Gamma :\mathcal{F}\rightarrow 2^{E}$ is an\textit{\
isotone operator } if for all $X,Y\in \mathcal{F}$ , $X\subseteq Y$ implies $%
\Gamma (X)\cap (E-Y)\subseteq \Gamma (Y)$. An accessible set system $(E,%
\mathcal{F})$ satisfies the interval property without upper bounds if and
only if $\Gamma :\mathcal{F}\rightarrow 2^{E}$ is an\textit{\ }isotone
operator: 
\[
x\in \Gamma (X)\cap (E-Y)\Leftrightarrow (x\in E-Y)\wedge (X\cup x\in 
\mathcal{F})\Rightarrow Y\cup x\in \mathcal{F}\Leftrightarrow x\in \Gamma
(Y). 
\]

\section{Isotone operators and antimatroids}

In this section a characterization of an antimatroid as a set system
generating by an isotone operator is given.

Consider an operator $\Psi :2^{E}\rightarrow 2^{E}$ such that 
\begin{equation}
for\ each\ X\subset E,\ \Psi (X)\subseteq E-X  \label{zero}
\end{equation}

In what follows we will use only operators satisfied (\ref{zero}). We can
build a set system, denoted by $(E,\mathcal{F}(\Psi ))$, using the following
algorithm.\newline

$\Psi $ \textbf{- Algorithm}

1. $\mathcal{F}(\Psi )=\emptyset $

2. $i=0$

3. $T_{i}=\{\emptyset \}$

4. while $T_{i}\neq \emptyset $

\qquad 4.1 $\mathcal{F}(\Psi )=\mathcal{F}(\Psi )\cup T_{i}$

$\qquad $4.2 $T_{i+1}=\{X\cup x$ : $(X\in T_{i})\wedge (x\in \Psi (X))\}$

\qquad 4.3 $i=i+1\bigskip $

Clearly, $(E,\mathcal{F}(\Psi ))$ is an accessible system, because for each
non empty $X\in \mathcal{F}(\Psi )$ there exists its ''parent'' $X-x\in 
\mathcal{F}(\Psi )$ by using which the set $X$ was generated.

It is not difficult to see that the property (\ref{zero}) implies that on
each step $i$ the algorithm generates the set $T_{i}$ in which each set has
exactly $i$ elements.

\smallskip Here are some examples of $\Psi $ operators:

$(a)$ Let for each $\ X\subset E,\ \Psi (X)=E-X$. Then $\mathcal{F}(\Psi
)=2^{E}$.

$(b)$ Let $E=\{x_{1},x_{2},...,x_{n}\}$, where for each $i$, $x_{i}<x_{i+1}$%
. Define for each $\emptyset \subset X\subset E,$ $\Psi (X)=\{x\in E:x>\max
(X)\}$, and $\Psi (\emptyset )=E$. It is easy to see that, $\mathcal{F}(\Psi
)=2^{E}$.

$(c)$ Let $(E,\leq )$ is a poset and $\Psi (X)=\{x\in E:x=\min (E-X)\}$,
then the obtained set system $(E,\mathcal{F}(\Psi ))$ is a poset antimatroid 
\cite{Greedoids}.

$(d)$ Let $E=\{x_{1},x_{2},...,x_{n}\}$, where for each $i$, $x_{i}<x_{i+1}$%
. Define for each $\emptyset \subset X\subset E,$ $\Psi (X)=\{x_{i+1}\}$,
where $x_{i}=\max (X)$, and $\Psi (\emptyset )=\{x_{1}\}$. Then $\mathcal{F}%
(\Psi )$ is a chain $\emptyset \subset \{x_{1}\}\subset
\{x_{1},x_{2}\}\subset ...\subset \{x_{1},x_{2},...,x_{n}\}$.

\smallskip Note, that the same set systems may be generated by different
operators (see the above examples). Now assume, that the operator $\Psi $ is
also an isotone operator, i.e.,

\begin{equation}
if\ X,Y\subset E\ then\ X\subseteq Y\ implies\ \Psi (X)\cap (E-Y)\subseteq
\Psi (Y).  \label{one}
\end{equation}

\begin{lemma}
\label{L-1}Let $\Psi $ be an isotone operator, then for each $X\in \mathcal{F%
}(\Psi )$ and $x\in E-X$, $X\cup x\in \mathcal{F}(\Psi )$ if and only if $%
x\in \Psi (X)$.
\end{lemma}

\smallskip \setlength {\parindent}{0.0cm}\textbf{Proof.} ''If'' immediately
follows from the structure of the $\Psi $-Algorithm. 
\setlength
{\parindent}{3.45ex}

Conversely, if set $X\cup x$ was generated from $X$ then $x\in \Psi (X)$. If
not - there is a sequence of sets generated by the $\Psi $-Algorithm $%
\emptyset =X_{0}\subset X_{1}\subset ...\subset X_{k}=X\cup x$ such that $%
X_{i}=X_{i-1}\cup x_{i}$ where $x_{i}\in \Psi (X_{i-1})$. Let $x=x_{j}$. Then

$(X_{j-1}\subseteq X)\wedge (x\in \Psi (X_{j-1}))\wedge (x\notin X)$ that
implies (from (\ref{one})) $x\in \Psi (X)$.\rule{2mm}{2mm}

\begin{corollary}
Two isotone operators $\Psi _{1}$ and $\Psi _{2}$ generate the same set
system $(E,\mathcal{F})$ if and only if $\Psi _{1}|_{\mathcal{F}}=\Psi
_{2}|_{\mathcal{F}}$.
\end{corollary}

The property (\ref{one}) makes possible to generate an antimatroid.

\begin{theorem}
\label{T-1}Set system $(E,\mathcal{F})$ is an antimatroid if and only if
there exists an isotone operator $\Psi $ such that $\mathcal{F}=\mathcal{F}%
(\Psi )$.
\end{theorem}

\setlength {\parindent}{0.0cm}\textbf{Proof.} Each $(E,\mathcal{F}(\Psi ))$
is an accessible system. If in additional $\Psi $ is an isotone operator,
then the set system $(E,\mathcal{F}(\Psi ))$ satisfies the interval property
without upper bounds. Indeed, if $X,Y\in \mathcal{F}(\Psi )$, and $%
X\subseteq Y$ and $x\in E-Y$, and $X\cup x\in \mathcal{F}(\Psi )$, then
(from Lemma \ref{L-1}) $x\in \Psi (X)$, that implies $x\in \Psi (Y)$, i.e. $%
Y\cup x\in \mathcal{F}(\Psi )$. It means (see Preposition \ref{P-1}) that
the set system $(E,\mathcal{F}(\Psi ))$ is an antimatroid. Moreover, $\Psi
(X)=\Gamma (X)$ for each $X\in \mathcal{F}(\Psi )$. 
\setlength
{\parindent}{3.45ex}

Conversely, let $(E,\mathcal{F})$ be an antimatroid. We will show that this
antimatroid can be generated by some isotone operator. First, build the
operator $\Psi :2^{E}\rightarrow 2^{E}$: 
\begin{equation}
for\ each\ X\subset E,\Psi (X)=\Gamma (\mathcal{B}_{X})  \label{two}
\end{equation}
where \-$\mathcal{B}_{X}$ is a basis of $X$. Since the basis is unique, the
definition is correct.

To show that the constructed operator is a required isotone operator we have
to prove two properties:

$(i)$ $\Psi (X)\subseteq E-X$

\quad Indeed, suppose that $x\in \Psi (X)$ and $x\in X$, then \-$\mathcal{B}%
_{X}\cup x\subseteq X$ and \-$\mathcal{B}_{X}\cup x\in \mathcal{F}$. Thus $%
\mathcal{B}_{X}\cup x$ is also a basis, a contradiction.

$(ii)$ $\Psi $ satisfies (\ref{one}).

\quad At first, if $X\subseteq Y$ then \-$\mathcal{B}_{X}\subseteq $\-$%
\mathcal{B}_{Y}$. Now, since $\Gamma $ is an isotone operator we have, by
using (\ref{two}), 
\[
(x\in E-Y)\wedge (x\in \Psi (X))\Rightarrow (x\in E-\mathcal{B}_{Y})\wedge
(x\in \Gamma (\mathcal{B}_{X}))\Rightarrow x\in \Gamma (\mathcal{B}%
_{Y})\Rightarrow x\in \Psi (Y) 
\]

It remains to show that $\mathcal{F}(\Psi )=\mathcal{F}$. For this purpose,
consider $X\in \mathcal{F}$. There exists a sequence $\emptyset
=X_{0}\subset X_{1}\subset ...\subset X_{k}=X$ where $X_{i}=X_{i-1}\cup
x_{i} $ and $X_{i}\in \mathcal{F}$ for $0\leq i\leq k$. Thus, $x_{i}\in \Psi
(X_{i-1})$, i.e. elements of the sequence is also obtained by $\Psi $%
-generator, and so $X\in \mathcal{F}(\Psi )$.

Conversely, let $X\in \mathcal{F}(\Psi )$. Then there is a sequence $%
\emptyset =X_{0}\subset X_{1}\subset ...\subset X_{m}=X$ where $%
X_{i}=X_{i-1}\cup x_{i}$ and $x_{i}\in \Psi (X_{i-1})$. Then $X_{i}\in 
\mathcal{F}$ for $0\leq i\leq m$,\- and so $X\in \mathcal{F}$.\rule{2mm}{2mm}

\section{The Chain Algorithm and monotone linkage functions}

In general, to optimize a set function is an \textit{NP}-hard problem, but
for some specific functions and for some specific set systems polynomial
algorithms are known. In this section we consider set functions defined as
minimum values of monotone linkage functions. Such set functions can be
maximized by a greedy type algorithm over a family of all subsets of $E$
(see \cite{Zaks}). Here we extend this result to antimatroids.

The monotone linkage functions were introduced by Mullat \cite{Mullat}. We
will give a necessary basic notions.

Let $\pi :E\times 2^{E}\rightarrow \mathbf{R}$ be a monotone linkage
function such that 
\begin{equation}
if\ X,Y\subseteq E\ and\ x\in E,\ then\ X\subseteq Y\ implies\ \pi (x,X)\geq
\pi (x,Y)  \label{three}
\end{equation}

For example, the single linkage $\pi (x,X)=\min_{y\in X}d_{xy}$, where $%
d_{xy}$ is a distance between two objects, is a monotone linkage function.

Consider $F:2^{E}\rightarrow \mathbf{R}$ defined for each $X\subset E$%
\begin{equation}
F(X)=\min_{x\in E-X}\pi (x,X)  \label{four}
\end{equation}

These functions were studied in \cite{Zaks},\cite{KMM}. A simple polynomial
algorithm which finds a set $X\subset E$ such that 
\[
F(X)=\max \{F(Y):Y\subset E\} 
\]
was developed, and the idea of this algorithm was used in searching of a
protein sequence alignment \cite{DIMACS}. In this section we extend our
results to a set system $(E,\mathcal{F}(\Psi ))$ generated by an isotone
operator $\Psi $. For this purpose we define a new set function as follows: 
\begin{equation}
F_{\Psi }(X)=\min_{x\in \Psi (X)}\pi (x,X)  \label{six}
\end{equation}

It should be pointed out that the definition (\ref{six}) is not limited to
set systems $(E,\mathcal{F}(\Psi ))$, but in order to the function $F_{\Psi
} $ to be defined on each subset $X\subset E$ the operator $\Psi $ must be
non-empty for each subset of $E$, i.e., 
\begin{equation}
for\ each\ X\subset E,\ \Psi (X)\neq \emptyset .  \label{seven}
\end{equation}

It is easy to show that a set system $(E,\mathcal{F}(\Psi ))$ with non-empty
isotone operator $\Psi $ is an antimatroid in which $E\in \mathcal{F}(\Psi )$%
. In \cite{Greedoids} this is a necessary condition for antimatroids,
whereas other authors doesn't involve this property in the definition of an
antimatroid. Thus, \cite{Bilbao} sets these antimatroids to the special
class of \textit{normal }antimatroids, and \cite{BZ} calls them \textit{full 
}antimatroids. In any case, an antimatroid $(E,\mathcal{F})$ has one and
only one maximal feasible set, namely $\cup_{X\in \mathcal{F}} X$, that we
denote as $E_{\mathcal{F}}$. Further we will only need the assumption that
operator $\Psi $ is not-empty on $\mathcal{F}-E_{\mathcal{F}} $.

\smallskip Consider the following optimization problem - given a monotone
linkage function $\pi $ , and a set system $(E,\mathcal{F}(\Psi ))$
generated by an isotone operator $\Psi $, find the feasible set $X\in 
\mathcal{F}(\Psi )-E_{\mathcal{F}(\Psi )}$ such that $F_{\Psi }(X)=\max
\{F_{\Psi }(Y):Y\in \mathcal{F}(\Psi )-E_{\mathcal{F}(\Psi )}\}$, where
function $F_{\Psi }$ defined by (\ref{six}). To solve this problem we build
the following algorithm.\newline

\textbf{The Chain Algorithm }$(E,\pi ,\Psi )$

1. Set $X^{0}=\emptyset $

2. Set $X=\emptyset $

3. While $\Psi (X)\neq \emptyset $ do

\qquad 3.1 If $F_{\Psi }(X)>F_{\Psi }(X^{0})$, set $X^{0}=X$

\qquad 3.2 Choose $x\in \Psi (X)$ such that $\pi (x,X)\leq \pi (y,X)$ for
all $y\in \Psi (X)$

\qquad 3.3 Set $X=X\cup x$

4. Return $X^{0}\bigskip $

Thus, the Chain Algorithm generates the chain of sets 
\[
\emptyset =X_{0}\subset X_{1}\subset ...\subset X_{k}=E_{\mathcal{F}(\Psi
)}, 
\]
where $X_{i}=X_{i-1}\cup x_{i}$ and $x_{i}\in \Psi (X_{i-1})$ for $1\leq
i\leq k$, and returns the minimal set $X^{0}$ of the chain on which the
value $F_{\Psi }(X^{0})$ is maximal.

\begin{theorem}
\smallskip \label{T-3}For a set system $(E,\mathcal{F}(\Psi ))$ the
following statements are equivalent

$(i)$ $\Psi $ is an isotone operator

$(ii)$ for all monotone linkage function $\pi $ the Chain Algorithm finds a
feasible set 

\quad that maximizes the function $F_{\Psi }$
\end{theorem}

\setlength {\parindent}{0.0cm}\textbf{Proof.} Let $X^{0}$ be the set
obtained by the Chain Algorithm. To prove that $X^{0}$ is a feasible set
maximizing $F_{\Psi }$, we have to show that $F_{\Psi }(X)\leq F_{\Psi
}(X^{0})$ for each $X\in \mathcal{F}(\Psi )-E_{\mathcal{F}(\Psi )}$. 
\setlength
{\parindent}{3.45ex}

Let $X_{0}\subset X_{1}\subset ...\subset X_{k}$ be the chain generated by
the Chain Algorithm. Let $j$ be the least integer for which $X_{j}\not%
{\subseteq }X$. Then $X_{j-1}\subseteq X$, $x_{j}\notin X$ and $X_{j-1}\cup
x_{j}\in \mathcal{F}(\Psi )$, that implies $x_{j}\in \Psi (X)$. Hence, 
\[
F_{\Psi }(X)\leq \pi (x_{j},X)\leq \pi (x_{j},X_{j-1})=F_{\Psi
}(X_{j-1})\leq F_{\Psi }(X^{0}). 
\]

Conversely, let $\Psi $ be not isotone operator, i.e. there exists $A,B\in 
\mathcal{F}(\Psi )-E_{\mathcal{F}(\Psi )}$ such that $A\subset B$, and there
is $a\in E-B$ such that $a\in \Psi (A)$ and $a\notin \Psi (B)$.
Accessibility of the set system $(E,\mathcal{F}(\Psi ))$ implies that there
exists a sequence 
\[
\emptyset =A_{0}\subset A_{1}\subset ...\subset A_{k}=A\subset A_{k+1}=A\cup
a, 
\]
where $A_{i}=A_{i-1}\cup a_{i}$ and $a_{i}\in \Psi (A_{i-1})$ for $1\leq
i\leq k$, and $a_{k+1}=a$. Define a monotone linkage function $\pi $ on
pairs $(x,X)$ where $X\subset E$ and $x\in E-X$: 
\[
\pi (x,X)=\left\{ 
\begin{array}{ll}
1, & X\supseteq A_{i-1}\ and\mathit{\ }x=a_{i}\ or\ A\cup a\subseteq
X\subset E\ and\ x\in E-X \\ 
2, & otherwise.
\end{array}
\right. 
\]
Then the Chain Algorithm generates a chain $A_{0}\subset ...\subset
A_{k}\subset A_{k+1}\subset ...\subset E_{\mathcal{F}(\Psi )}$, on which the
values of the function $F_{\Psi }$ are equal to $1$, but $F_{\Psi }(B)=2$.
Thus, the Chain Algorithm does not find a feasible set that maximizes the
function $F_{\Psi }$. \rule{2mm}{2mm}

The Chain Algorithm is a greedy type algorithm since it based on the best
choice principle: it chooses on each step the extreme element (in sense of
linkage function) and thus approaches the optimal solution. Let $P$ is the
maximum complexity of $\pi (x,X)$ computation over all pairs $(x,X)$ where $%
x\in E-X$. Then the Chain Algorithm finds the optimal feasible set in $%
O(P|E|^{2})$ time, for example, in some clustering problems \cite{KMM}, the
complexity of the Chain Algorithm is $O(|E|^{3})$.

Notice, that for antimatroids the functions $\Psi $ and $\Gamma $ are
identical (see Theorem \ref{T-1}) then the following central result
immediately follows from previous theorems:

\begin{theorem}
For an accessible set system $(E,\mathcal{F})$, where $\Gamma (X)\neq
\emptyset $  for each  $X\in \mathcal{F}-E_{\mathcal{F}}$ the following
statements are equivalent

$(1)$ the set system $(E,\mathcal{F})$ is an antimatroid

$(2)$ The Chain Algorithm finds a feasible set that maximizes the function $%
F_{\Gamma }$ for 

\quad every monotone linkage function $\pi $
\end{theorem}

\section{Truncated antimatroids}

In this section we extend our results obtained in Section 4 to \textit{%
truncated antimatroids} considered in the work of Boyd and Faigle \cite{BF}.

\begin{definition}
The \textit{k-truncation }of a set system $(E,\mathcal{F})$ is a set system
defined by 
\[
\mathcal{F}_{k}=\{X\in \mathcal{F}:\left| X\right| \leq k\}
\]
\end{definition}

If $\mathcal{F}$ is an antimatroid, then $\mathcal{F}_{k}$ is a k-truncated
antimatroid.

The \textit{rank} of a set $X\subseteq E$ is defined as $\varrho (X)=\max
\{\left| Y\right| :(Y\in \mathcal{F})\wedge (Y\subseteq X)\}$, the rank of
set system $(E,\mathcal{F})$ is defined as $\varrho (\mathcal{F})=\varrho
(E).$ For a given antimatroid $(E,\mathcal{F})$ the rank of k-truncated
antimatroid $\varrho (\mathcal{F}_{k})=k$, whenever $k\leq \varrho (\mathcal{%
F})$.

Let $(E,\mathcal{F}(\Psi ))$ be an antimatroid generated by an isotone
operator $\Psi $. Consider a (\textit{k-1)-truncated} operator 
\begin{equation}
\Psi _{k-1}(X)=\left\{ 
\begin{array}{ll}
\Psi (X) , & \left| X\right| \leq k-1 \\ 
\emptyset , & otherwise
\end{array}
\right.  \label{def_oper}
\end{equation}

The set system generated by $\Psi _{k-1}$-operator is a $k$-truncated
antimatroid$\ (E,(\mathcal{F}(\Psi ))_{k})$, i.e., $(\mathcal{F}(\Psi ))_{k}=%
\mathcal{F}(\Psi _{k-1})$. Indeed, any set $\left| Y\right| \leq k$ is
generated by $\Psi _{k-1}$-generator if and only if $Y\in \mathcal{F}$ ,
since $\Psi _{k-1}(X)\equiv \Psi (X)$ for all $\left| X\right| \leq k-1$.
Moreover, assume a set $Y$ for which $\left| Y\right| >k$ was obtained by $%
\Psi _{k-1}$-generator, then there is a set $X$ such that $\left| X\right|
\geq k$ and $\Psi _{k-1}(X)\neq \emptyset $, in contradiction with the
definition (\ref{def_oper}).

Clearly, that the $\Psi _{k-1}$-operator is not isotone on all $2^{E}-E$,
but it satisfies to the following condition: 
\begin{equation}
X,Y\subset E,\ (X\subseteq Y)\wedge (\left| Y\right| \leq k-1)\ implies\
\Psi _{k-1}(X)\cap (E-Y)\subseteq \Psi _{k-1}(Y)  \label{trunc_isot}
\end{equation}

\smallskip We call an operator $\Psi $ a\textit{\ (k-1)-isotone operator }if
it satisfies (\ref{trunc_isot}) and $\Psi (X)=\emptyset $ for each $X\subset
E$, such that $\left| X\right| \geq k$.

The following theorem shows that a (k-1)-isotone operator generates a
truncated antimatroids in the same way as an isotone operator determines an
antimatroid:

\begin{theorem}
Set system $(E,\mathcal{F})$ is a k-truncated antimatroid if and only if
there exists a (k-1)-isotone operator $\Psi $ such that $\mathcal{F}=%
\mathcal{F}(\Psi )$.
\end{theorem}

\setlength {\parindent}{0.0cm}\textbf{Proof.} To prove that the set system $%
(E,\mathcal{F}(\Psi ))$ is a k-truncated antimatroid we have to build an
antimatroid of which it is a truncation. Define, by analogy with \cite{BF} 
\begin{equation}
\Omega =\{X\subseteq E: there\ are\ some\ X_{1},...,X_{p}\in \mathcal{F}%
(\Psi ) such\ that\ X=X_{1}\cup ...\cup X_{p}\},  \label{closeU}
\end{equation}
i.e., $\Omega $ is a closure by union of $\mathcal{F}(\Psi )$. 
\setlength
{\parindent}{3.45ex}

The set system $(E,\Omega )$ is closed under union, so to prove that it is
an antimatroid it is remain to verify that the set system $(E,\Omega )$ is
accessible. By analogy with \cite{BF} consider a set $X\in \Omega $ and let $%
X=X_{1}\cup ...\cup X_{k}$. Then there exists $x\in X_{1}$ such that $%
X_{1}-x\in \mathcal{F}(\Psi )$. Assume without loss of generality that $%
x\notin X_{2},X_{3},...,X_{k}$, for otherwise we could let $X_{1}=X_{1}-x$.
If so, $X-x=(X_{1}-x)\cup X_{2}\cup ...\cup X_{k}\in \Omega $.

To show that the k-truncation of $(E,\Omega )$ is $(E,\mathcal{F}(\Psi ))$
it is sufficient to prove that $X\in \mathcal{F}(\Psi )$ if and only if $%
X\in \Omega $ and $\left| X\right| \leq k$. Indeed, if $X\in \mathcal{F}%
(\Psi )$ then, from (\ref{closeU}), $X\in \Omega $ and obviously $\left|
X\right| \leq k$. Conversely, let $X\in \Omega $ and $\left| X\right| \leq k$%
, then $X=A_{1}\cup ...\cup A_{p}$. We show that $X\in \mathcal{F}(\Psi )$
by induction on $p$. If $p=1$ then $X\in \mathcal{F}(\Psi )$. Consider $%
A=A_{1}\cup ...\cup A_{p-1}$. By the hypothesis of induction, $A\in \mathcal{%
F}(\Psi )$ . Assume $|A|<k$, for otherwise $X=A$ and then $X\in \mathcal{F}%
(\Psi )$. Let $\emptyset =X_{0}\subset X_{1}\subset ...\subset X_{l}=A$ be a
sequence of sets generated by the $\Psi $-operator, where $X_{i}=X_{i-1}\cup
x_{i}$ and $x_{i}\in \Psi (X_{i-1})$ for $1\leq i\leq l<k$. Let $j$ be the
least integer for which $X_{j}\not{\subseteq}A_{p}$. Then $X_{j-1}\subseteq
A_{p}$, $x_{j}\notin A_{p}$ and $X_{j-1}\cup x_{j}\in \mathcal{F}(\Psi )$ ,
that implies $x_{j}\in \Psi (A_{p})$. Repeated application of (\ref
{trunc_isot}) yields a set $X=A_{p}\cup (A-A_{p})\in \mathcal{F}(\Psi )$.

Conversely, let $(E,\mathcal{F})$ be a k-truncated antimatroid, then there
is an antimatroid $(E,\Phi )$ for which $(E,\mathcal{F})$ is a k-truncation.
Since $(E,\Phi )$ is an antimatroid, there exists (Theorem \ref{T-1}) an
isotone operator $\widehat{\Psi }$ that generates the antimatroid $(E,\Phi )$%
, and (k-1)-truncation $\widehat{\Psi _{k-1}}$ generates the set system $(E,%
\mathcal{F})$. Obviously, the operator $\widehat{\Psi _{k-1}}$ satisfies to (%
\ref{trunc_isot}). \rule{2mm}{2mm}

Now, using the same technique as in Theorem \ref{T-3} we obtain an
algorithmic characterization of truncated antimatroids.

\begin{theorem}
For an accessible set system $(E,\mathcal{F})$, where $\Gamma (X)\neq
\emptyset $ for each $X\in \mathcal{F}$ such that $|X|<k$ the following
statements are equivalent

$(1)$ the set system $(E,\mathcal{F})_{k}$ is a k-truncated antimatroid

$(2)$ the Chain Algorithm finds a feasible set that maximizes the function $%
F_{\Gamma }$ on 

$\quad (E,\mathcal{F})_{k-1}$ for any monotone linkage function $\pi $
\end{theorem}

\section{Correspondence between two algorithmic characterization of
antimatroids}

In this section we consider another algorithmic approach to antimatroids
introduced in work of Boyd and Faigle \cite{BF}. Since the work based on
other definition of an antimatroid as a formal language, some additional
notation is needed. Given a finite alphabet $E$ consists of \textit{letters}%
. A \textit{word }over $E$ is a sequence of letters from $E$, denoted by the
lower case of Greek letters $\alpha $,$\beta $ and $\gamma $. A \textit{%
language }$\mathcal{L}$ is a set of words of $E$. The concatenation of two
words $\alpha $ and $\beta $ will be denoted $\alpha \beta $, $\alpha _{k}$
will be used to denote a word of length $k$ and the set of distinct letters
in a word $\alpha $ will be denoted $\widetilde{\alpha }$. The language is
called \textit{simple} if there are no words with repeated letters.

\begin{definition}
An antimatroid language is a simple language $(E,\mathcal{L})$ satisfying
the following two properties:

$(1)$ If $\alpha x\in \mathcal{L}$, then $\alpha \in \mathcal{L}$.

$(2)$ If $\alpha $,$\beta \in \mathcal{L}$ and $\widetilde{\alpha }\not%
{\subseteq}\widetilde{\beta }$, then there exists an $x\in \widetilde{\alpha 
}$ such that $\beta x\in \mathcal{L}$.
\end{definition}

Antimatroids and antimatroid languages are equivalent in the following sense 
\cite{Greedoids}.

\begin{theorem}
If $(E,\mathcal{L})$ is an antimatroid language, then 
\[
F(\mathcal{L})=\{\widetilde{\alpha }:\alpha \in \mathcal{L}\} 
\]

is an antimatroid $(E,F(\mathcal{L}))$.

Conversely, if $(E,\mathcal{F})$ is an antimatroid, then 
\[
L(\mathcal{F})=\{x_{1}...x_{k}:\{x_{1},...x_{j}\}\in \mathcal{F\ }for\ 1\leq
j\leq k\} 
\]

is an antimatroid language $(E,L(\mathcal{F}))$. Further, $L(F(\mathcal{L}))=%
\mathcal{L}$ and $F(L(\mathcal{F}))=\mathcal{F}$.
\end{theorem}

The next problem was considered in \cite{BF}: let $f:E\times
2^{E}\rightarrow \mathbf{R}$ be a monotone function such that $f(x,A)\leq
f(x,B)$ whenever $B\subseteq A$. Define a \textit{maximum nesting function } 
\[
W(x_{1}...x_{k})=\max
\{f(x_{1},\{x_{1}\}),...,f(x_{k},\{x_{1},...,x_{k}\})\}. 
\]

\textit{The minimax nesting problem} was defined as follows: given a simple
language $(E,\mathcal{L})$ with a monotone function $f$ and a nonnegative
integer $k\leq \varrho (\mathcal{L})$, find $\alpha _{k}\in \mathcal{L}$
such that 
\[
W(\alpha _{k})=\min \{W(\beta _{k}):\beta _{k}\in \mathcal{L}\} 
\]

The main theorem proved in \cite{BF} is reads as follows.

\begin{theorem}
Let $(E,\mathcal{L})$ be a simple language. The greedy algorithm solves the
minimax nesting problem for every monotone function $f$ if and only if $(E,%
\mathcal{L})$ is a truncated antimatroid.
\end{theorem}

We will show the correspondence between our algorithmic characterization of
antimatroids and characterization of Boyd and Faigle.

First note, that in \cite{BF} was proved that the constructed word $\alpha
_{k}=x_{1}...x_{k}$ satisfies also the following property: 
\begin{equation}
W(x_{1}...x_{i})=\min \{W(\beta _{i}):\beta _{i}\in \mathcal{L}\}\ for\
each\ i\ such\ that\ 1\leq i\leq k  \label{last}
\end{equation}

Second note, that the Chain Algorithm builds a sequence $\emptyset
=X_{0}\subset X_{1}\subset ...\subset X_{k}$ where $X_{i}=X_{i-1}\cup x_{i}$
for $1\leq i\leq k$, i.e. this algorithm generates the sequence $%
x_{1}...x_{k}$. So all sets $X_{i}$, obtained by the Chain Algorithm, has a
natural order: $X_{i}=\{x_{1},...,x_{i}\}$, i.e. we can consider each set $%
X_{i}$ also as a word $\alpha _{i}=x_{1}...x_{i}$. We are now ready to prove:

\begin{theorem}
\label{T_last}Let $(E,\mathcal{L})$ be a k-truncated antimatroid and $\Psi $
is a operator which generates this antimatroid , let 
\[
f(x_{i},\{x_{1},...,x_{i}\})=\pi (x_{i},\{x_{1},...,x_{i-1}\})\ for\ each\
i\ such\ that\ 1\leq i\leq k\, 
\]

then

$(i)$ if $X^{0}$ is an optimal set obtained by the Chain Algorithm, then
there exists a word $\alpha _{k}\in \mathcal{L}$ that satisfies (\ref{last})
and $X^{0}=\{x_{1},...,x_{p}\}$ is a shortest prefix of $\alpha _{k}$ such
that $W(x_{1}...x_{p+1})=W(\alpha _{k})=F_{\Psi }(X^{0})$.

$(ii)$ if $\alpha _{k}$ is a solution of the minimax nesting problem
obtained by the greedy algorithm, then a shortest prefix $%
\{x_{1},...,x_{p}\} $ of $\alpha _{k}$ such that $W(x_{1}...x_{p+1})=W(%
\alpha _{k})$ maximizes the function $F_{\Psi }$.
\end{theorem}

\setlength {\parindent}{0.0cm}\textbf{Proof.} $(i)$ Let $x_{1}...x_{k}$ be
the sequence generating by the Chain Algorithm and let $X^{0}=%
\{x_{1},...,x_{p}\}.$ Set $\alpha _{k}=x_{1}...x_{k}$ and prove that $\alpha
_{k}$ satisfies (\ref{last}). Suppose not, then let $\gamma
_{m}=y_{1}...y_{m}$ be a shortest word such that $W(\gamma
_{m})<W(x_{1}...x_{m})$. It means that for each $i<m$%
\[
\max \{\pi (x_{1},\emptyset ),...,\pi (x_{i},\{x_{1},...,x_{i-1}\})\}\leq
\max \{\pi (y_{1},\emptyset ),...,\pi (y_{i},\{y_{1},...,y_{i-1}\})\} 
\]

and for each $i\leq m$%
\begin{equation}
\pi (x_{m},\{x_{1},...,x_{m-1}\})>\max \{\pi (y_{1},\emptyset ),...,\pi
(y_{i},\{y_{1},...,y_{i-1}\})\}  \label{last2}
\end{equation}
\setlength
{\parindent}{3.45ex}

If $\{y_{1}...y_{m-1}\}=\{x_{1},...,x_{m-1}\}$, then $y_{m}\in \Psi
(\{x_{1},...,x_{m-1}\})$ and from (\ref{last2}) 
\[
\pi (y_{m,}\{x_{1},...,x_{m-1}\})=\pi (y_{m},\{y_{1},...,y_{m-1}\})<\pi
(x_{m},\{x_{1},...,x_{m-1}\}) 
\]
So the Chain Algorithm should choose $y_{m}$ and not $x_{m}$.

Thus, let $j$ be the smallest index such that $\{y_{1},...,y_{j-1}\}%
\subseteq $ $\{x_{1},...,x_{m-1}\}$ and $y_{j}\notin \{x_{1},...,x_{m-1}\}$.
Since $y_{j}\in \Psi (\{y_{1},...,y_{j-1}\})$, we get that $y_{j}\in \Psi
(\{x_{1},...,x_{m-1}\})$. Hence, from monotonicity of $\pi $ and from (\ref
{last2}) 
\[
\pi (y_{j},\{x_{1},...,x_{m-1}\})\leq \pi (y_{j},\{y_{1},...,y_{j-1}\})<\pi
(x_{m},\{x_{1},...,x_{m-1}\}) 
\]
contradiction to optimal choice of $x_{m}$.

Finally, the Chain Algorithm construction implies, that $X^{0}=%
\{x_{1},...,x_{p}\}$ is the shortest prefix of $\alpha _{k}$ such that 
\[
F_{\Psi }(X^{0})=\pi (x_{p+1},\{x_{1}...x_{p}\})=W(x_{1}...x_{p+1})=W(\alpha
_{k}) 
\]

$(ii)$ Conversely, let $\alpha _{k}$ be a solution of the minimax nesting
problem and let $X^{0}=x_{1},...,x_{p}$ be the shortest prefix such that $%
W(x_{1}...x_{p+1})=W(\alpha _{k})$. Then 
\[
\pi (x_{p+1},\{x_{1}...x_{p}\})>\pi (x_{i+1},\{x_{1}...x_{i}\})\ for\ i<p 
\]
and 
\[
\pi (x_{p+1},\{x_{1}...x_{p}\})\geq \pi (x_{i+1},\{x_{1}...x_{i}\})\ for\
i\geq p 
\]

Certainly, $\pi (x_{p+1},\{x_{1}...x_{p}\})=\min_{x\in \Psi (X^{0})}\pi
(x,\{x_{1}...x_{p}\})$. If not, there is $x^{0}\in \Psi (X^{0})$ such that $%
\pi (x^{0},\{x_{1}...x_{p}\})<\pi (x_{p+1},\{x_{1}...x_{p}\})$, i.e. $%
W(x_{1}...x_{p}x^{0})<W(x_{1}...x_{p+1})$ - contradiction to (\ref{last}).
So, $F_{\Psi }(X^{0})=\pi (x_{p+1},\{x_{1}...x_{p}\})$.

Consider some set $X\in F(\mathcal{L})$. If $X=\{x_{1}...x_{j}\}$, i.e. $X$
is some prefix of $\alpha _{k}$, then 
\[
F_{\Psi }(X)=\min_{x\in \Psi (X)}\pi (x,X)\leq \pi
(x_{j+1},\{x_{1}...x_{j}\})\leq \pi (x_{p+1},\{x_{1}...x_{p}\})=F_{\Psi
}(X^{0}) 
\]

Otherwise, let $j$ be the smallest index such that $\{x_{1}...x_{j}\}%
\subseteq X$ and $x_{j+1}\notin X$. Then $x_{j+1}\in \Psi (X)$. Hence, 
\[
F_{\Psi }(X)=\min_{x\in \Psi (X)}\pi (x,X)\leq \pi (x_{j+1},X)\leq 
\]
\[
\leq \pi (x_{j+1},\{x_{1}...x_{j}\})\leq \pi
(x_{p+1},\{x_{1}...x_{p}\})=F_{\Psi }(X^{0}). 
\]

\rule{2mm}{2mm}

\section{Conclusions}

In this article, we discussed a set system algorithmic description of one
subclass of greedoids, namely, antimatroids. Further we compared a new
description with a known one based on the approach to define greedoids as
languages. Actually, there are some more important subclasses of greedoids
also enjoying natural algorithmic characterizations in terms of their
feasible set systems, for instance, matroids and Gaussian greedoids. These
results may lead to new algorithmic frameworks for additional types of
greedoids. We consider the family of interval greedoids as a strong
candidate for the collection of successes of the set system algorithmic
approach.

\end{document}